\begin{document} 
\selectlanguage{english}
\title{Calabi-Yau threefolds in positive characteristic}
\author{Yukihide Takayama}
\affilation{
	Prof. Dr. Yukihide Takayama \\
	Department of Mathematical Sciences, Ritsumeikan University\\
        1-1-1 Nojihigashi, Kusatsu, Shiga, 525-8577, Japan\\
	\url{takayama@se.ritsumei.ac.jp} 
}
\smallFoot
\maketitle
{\it 
  Expanded lecture notes of a course given as part of
  XVI International Workshop for Young Mathematicians
  ``Algebraic Geometry''
  held in Krakow,  September 18 -- 24, 2016.}

\begin{abstract}
  In this note, an overview of Calabi-Yau varieties in
  positive characteristic is presented. Although Calabi-Yau varieties in
  characteristic zero are unobstructed, there are examples of
  Calabi-Yau threefolds in positive characteristic which cannot be
  lifted to characteristic zero, although one-dimensional and
  two-dimensional Calabi-Yau varieties, i.e., elliptic curves and K3
  surfaces, are all liftable to characteristic zero. In this respect,
  Calabi-Yau threefolds in positive characteristic are interesting in
  view of deformation theory and they are still very mysterious.\\
\end{abstract}



%
%
\usepackage{amscd}
\def\opn#1#2{\def#1{\operatorname{#2}}} 
\opn\chara{char}
\opn\Spf{Spf}
\opn\Spec{Spec}
\opn\Ker{Ker}
\opn\Im{Im}
\opn\Hom{Hom}
\opn\Ext{Ext}
\opn\rk{rk}
\opn\Sing{Sing}
\def\mm{{\frk m}}
\def\OO{\mathcal{O}}
\def\cal#1{\mathcal{#1}}
\let\iso=\cong
\let\Dirsum=\bigoplus
\let\tensor=\otimes
\let\dirsum=\oplus
\let\To=\longrightarrow
\def\frk{\frak}               
\def\NZQ{\Bbb}               
\def\QQ{{\NZQ Q}}
\def\ZZ{{\NZQ Z}}
\def\PP{{\NZQ P}}
\def\CC{{\NZQ C}}
\def\FF{{\NZQ F}}

\tableofcontents

\newpage
\section{Introduction}

Calabi-Yau varieties, in particular Calabi-Yau 3-folds, over complex
numbers have been extensively studied over the last decades mostly from the
interest in mathematical physics. The most notable
property of complex Calabi-Yau varieties is that they are unobstructed
in deformation.

However, for Calabi-Yau varieties of dimension $\geq 3$ in positive
characteristic, the situation is quite different. In late 1990's,
M.~Hirokado found an example of Calabi-Yau 3-fold in characteristic $3$
that cannot be lifted to characteristic $0$ \cite{H99}. After that
other examples of non-liftable Calabi-Yau 3-fods have been found
\cite{H01, Schr04, HISark,HISmanuscripta}.  They are all in
characteristic $p=2$ and $p=3$. However, D.~van Straten, S.Cynk and
M.~Sch\"utt \cite{CyvS, CySchu} found a large number of examples of
non-liftable Calabi-Yau algebraic spaces in characteristic $p\geq 5$,
which are not schemes anymore.  It is still an open problem whether
there exist non-liftable Calabi-Yau varieties in characteristic $\geq
5$.

On the other hand, by the cerebrated theorem by P.~Deligne and
L.~Illusie (and M.~Raynaud) \cite{DI87}, for a projective variety $X$
over an algebraically closed field $k$ of $\chara(k)=p\geq \dim X$, if 
$X$ can be lifted to the ring $W_2(k)$ of second Witt vectors,
Hodge-to-de Rham spectral sequence degenerates at $E_1$ and
Kodaira-Akizuki-Nakano vanishing of cohomologies holds. We note that
$W_2$-liftability is not a necessary condition for Kodaira type
vanishing and there are examples of varieties that are $W_2$-non-liftable
but Kodaira vanishing holds.

Thus, for non-liftable Calabi-Yau varieties,
even if they are not liftable over the ring $W(k)$ of Witt vectors,
which implies liftability to characteristic $0$, it is an interesting
question whether they are liftable over $W_2(k)$.
The above mentioned Hirokado's  example and
the Schr\"oer's examples  are known to be non-liftable over
$W_2(k)$ \cite{Ek03} and $W_2$-liftability is still open for other
examples.  So Kodaira type vanishing for non-liftable Calabi-Yau varieties in
positive characteristic is still a mysterious problem.

In this note, we give an overview of the research on Calabi-Yau
3-folds in positive characteristic in the last decades.  The
configuration of this note is as follows.  In section~2, we
summarize briefly what is different from the geometry of Calabi-Yau varieties in
characteristic~$0$.  Among the unique features of geometry in
positiv characteristic, we will focus on obstructedness of
deformation, which we will elaborate in section~3. Final section is
devoted to construction of non-liftable Calabi-Yau 3-folds. We first
summarize the known reasons for non-liftability and overview the
examples that have been found so far. We also mention what is known
about Kodaira type vanishing for non-liftable Calabi-Yau varieties.

\section{Unique features of CY 3-folds in positiv characteristic}

In this section, we overview unique features of Calabi-Yau
3-folds in positive characteristic as compared to
characteristic $0$ case. Some of the feature will be
considered in detail in the subsequent sections.

\begin{definition}[Calabi-Yau $n$-fold]
  A {\em Calabi-Yau $n$-fold} $X$ is a smooth projective 
  variety over an algebraically closed field $k$ such that
  $H^i(X,\OO_X)=0$ for $0<i<n=\dim{X}$ and the canonical sheaf is
  trivial $\omega_X \iso \OO_X$.  In particular, a Calabi-Yau $1$-fold
  is an elliptic curve and a Calabi-Yau $2$-fold is a K3 surface.
\end{definition}

In the following, we will mostly consider the case of $n=3$ and
$\chara(k)=p>0$.
Some authors assume only properness instead of projectivity for
Calabi-Yau varieties. But in this note we will always assume
projectivity for a Calabi-Yau variety.

\subsection{Hodge diamond without Hodge symmetry}

By Serre duality, we know that 
the Hodge numbers $h^{ij} := \dim_k H^k(X,\Omega_X^i)$
of a Calabi-Yau 3-fold $X$ are as follows:
\begin{equation*}
  \begin{array}{ccccccc}
      &       &        &  h^{00} &        &       &       \\
      &       & h^{10}  &        & h^{01}  &       &       \\
      & h^{20} &        & h^{11}  &        & h^{02} &       \\
h^{30} &       & h^{21}  &        & h^{12}  &       & h^{03} \\
      & h^{31} &        & h^{22}  &        & h^{13} &        \\
      &       & h^{32}  &        & h^{23}  &       &       \\
      &       &        & h^{33}  &        &       &       \\
  \end{array}
\qquad = \qquad 
    \begin{array}{ccccccc}
      &       &        &  1     &        &       &       \\
      &       & h^{10}  &        & 0      &       &       \\
      & h^{20} &        & h^{11}  &        & 0     &       \\
      1       &        & h^{12}  &        & h^{12}  &       & 1 \\
      & 0     &        & h^{11}  &        & h^{20} &        \\
      &       & 0      &        & h^{10}  &       &       \\
      &       &        & 1      &        &       &       \\
  \end{array}
  \end{equation*}
In positive characteristic, Hodge symmetry ($h^{ij}=h^{ji}$)
does not hold in general and we have 
\begin{proposition}[Hodge symmetry]
  \label{hodgesymm}
  For a Calabi-Yau 3-fold $X$, Hodge symmetry holds if and only if
  $h^{10}=h^{20}=0$, namely $H^0(X,\Omega_X^1)=H^0(X, T_X)=0$, where
  $T_X = {\cal Hom}_{\OO_X}(\Omega_X^1, \OO_X)$ is the tangent bundle.
\end{proposition}
\begin{proof}
The last part uses the isomorphism
 $\Omega_X^2 \iso \Omega_X^3 \wedge \Omega_X^{-1}
= \OO_X \wedge \Omega_X^{-1}\iso T_X$.
\end{proof}

\subsection{obstructed deformation}

\subsubsection{characteristic $0$ case}
We first recall the notion of universal deformation.

\begin{definition}[deformation]
  Let $X$ be a complex manifold.
  \begin{enumerate}
    \item A {\em deformation} of $X$ is
  a smoot proper morphism $\psi: {\cal X}\to (S,0)$, where
  ${\cal X}$ and $S$ are connected complex spaces and
  $0\in S$ a distinguished point, 
  such that $X\iso {\cal X}_0:= \psi^{-1}(0)$, 
\item
  A deformation ${\mathcal X}\to (S,0)$ is
  called {\em universal} if any other deformation
  $\psi':{\mathcal X}'\to (S',0')$ is isomorphic to
  the pull-back under a uniquely determined morphism
  $\varphi:S'\To S$ with $\varphi(0')=0$.
  We denote a universal deformation
  by ${\mathcal X} \To \operatorname{Def}(X)$
  \end{enumerate}
 \end{definition}

\begin{theorem}[Bogomolov-Tian-Todorov]
  \label{BTT}
  Let $X$ be a Calabi-Yau manifold of any dimension over algebraically
  closed field $k$ of $\chara(k)=0$.  Then $\operatorname{Def}(X)$ is a germ of a
  smooth manifold with tangent space $H^1(X, T_X)$.
  \end{theorem}
\begin{proof}
  By Lefschetz principle and GAGA, we may consider $X$ as a
  compact complex K\"ahler manifold.  Then we prove smoothness of
  $\operatorname{Def}(X)$
  by a complex analytic method. For the detail, see Theorem~14.10~\cite{GHJ}
  and the references cited there.
\end{proof}

According to deformation theory, obstruction is contained in $H^2(X,
T_X)$. So if $H^2(X, T_X)=0$ we can say that $X$ is unobstructed.
But Theorem~\ref{BTT} claims that even if $H^2(X,T_X)\ne{0}$,
which is actually possible in dimension $\geq 3$, 
its elements are not an  obstruction to deformation.

\begin{remark}
In the algebraic setting,
Theorem~\ref{BTT} means that any deformation of $X$ is unobstructed
in the sense that for any small extension $\varphi: B\To A$, namely,
for any surjective homomorphism
$\varphi$ of local Artinian $\CC$-algebras with $\mm_B\Ker\varphi=0$,
any variety $X_A$ over $\Spec(A)$ can be lifted over $\Spec{B}$.
\end{remark}

\subsubsection{Elliptic curves}

For curves, we have $H^2(X, T_X)=0$ for the dimensional reason
so that smooth projective curves in any characteristic are unobstructed.
In particular elliptic curves are unobstructed.

\subsubsection{K3 surfaces}

For a K3 surface $X$, we have $H^2(X, T_X)=0$ by the following theorem,
so that it is unobstructed.

\begin{theorem}[Deligne~\cite{DI81}]
  \label{DI:K3}
  Let $X$ be a K3 surface over an algebraically closed field $k$
  of $\chara(k)=p>0$. Then
  \begin{enumerate}
  \item Hodge to de Rham spectral sequence
    \begin{equation*}
      E_1^{pq} = H^q(X, \Omega_X^p)\Rightarrow H_{DR}^{p+q}(X/k)
      \end{equation*}
    degenerates at $E_1$. In particular,
    we have the Hodge decomposition
    \begin{equation*}
      H^k_{DR}(X/k) \iso \sum_{i=0}^k H^{k-i}(X, \Omega_X^i)
    \end{equation*}
  \item the Hodge diamond is
    \begin{equation*}
            \begin{array}{ccccc}
          &       & h^{00} &       & \\
          & h^{10} &       & h^{01} &   \\
    h^{20} &       & h^{11}&        & h^{02} \\
          & 0     &       & h^{12} &   \\
          &       & h^{22} &   & 
            \end{array}
\quad  = \quad
      \begin{array}{ccccc}
          &   & 1  &   & \\
          & 0 &    & 0 &   \\
        1 &   & 20 &   & 1 \\
          & 0 &    & 0 &   \\
          &   & 1  &   & 
        \end{array}
    \end{equation*}
    \end{enumerate}
\end{theorem}

Unobstructedness $H^2(X,T_X)=0$ in Theorem~\ref{DI:K3} implies
that a K3 surface over $k$ can be formally lifted over the ring $W(k)$
of Witt vectors, i.e.,
there exists a formal scheme $\hat{{\mathcal X}}$ over $W(k)$ such
that $\hat{{\mathcal X}}\times_{W(k)}k \iso X$. Moreover,
$\dim_kH^1(X,T_X)=20$ means that the moduli space of formal K3
surfaces is 20 dimensional.

On the other hand, if we want an {\em algebraic} lifting of a K3
surface, the situation is a little subtler. 

\begin{quote}
\underline{{\bf Projective lifting problem:}} Let $X$ be a projective 
variety over a field $k$ of positive characteristic.
Then, find a projective scheme ${\mathcal X}$
over $S=\Spec R$, where $R$ is a ring of mixed characteristic,
such that $X \iso {\mathcal X}\times_S k$.
\end{quote}
Recall that a ring $R$ of mixed characteristic means an integral
domain in $\chara(R)=0$ with a maximal ideal $\mm \subset R$
satisfying $\chara \left(R/\mm\right) > 0$.  Here, we can consider
$R=W(k)$, for example. In this case, we also need to lift an ample
line bundle on $X$ to the formal lifting $\hat{{\mathcal X}}$ and apply
Grothendieck's algebraization theorem
(see Th\'eor$\grave{{\rm e}}$me~(5.4.5) \cite{EGAIII}
or Theorem~21.2~\cite{HDT}) to obtain a projective scheme $\tilde{X}$ whose
formal completion at the closed fiber $X$ is $\hat{{\mathcal X}}$.

The obstruction to lifting an invertible sheaf is contained  in
$H^2(X,\OO_X)$ (see Theorem~\ref{HDT:th22.1} below)
which is $1$-dimensional for a K3 surface.
Thus, the moduli space of {\em algebraic} K3 surfaces is
smaller by one dimension, namely $19$ dimensional.

Projective lifting problem for K3 surfaces 
has been solved except the case $p=2$.

\begin{theorem}[Ogus~\cite{Og79}]
For $p>2$, a K3 surface $X$ can be projectively lifted over $W(k)$.
\end{theorem}
\begin{proof}
  By Corollary~2.3~\cite{Og79}, a K3 surface $X$ can be lifted over
  $W(k)$ if $X$ is not ``superspecial''. By Remark~2.3~\cite{Og79}, if
  Tate conjecture for smooth proper surfaces \cite{Ta65} holds,
  the only ``superspecial'' K3 surface is the
  Kummer surface associated to a product of supersingular elliptic
  curves and we can show that this can be lifted over $W(k)$.
  Finally, Tate's conjecture has been established for $p\geq 3$
  by several authors.
 \end{proof}

See, for example \cite{Lied16}, for the detail of deformation theory
of K3 surfaces in positive characteristic.

\subsubsection{CY $n$-folds ($n\geq 3$)}

Apart from the cases of $\chara(k)=0$ or $\chara(k)=p>0$ with $\dim
\leq 2$, which we have seen so far, the situation is
quite different for $\dim \geq 3$. Namely,

\begin{theorem}[Hirokado, Schr\"oer, Ito, Saito, Ekedahl,
    Cynk, van Straten, Sch\"utt]
  There are examples of Calabi-Yau 3-folds over
  an algebraically closed field $k$ of characteristic $p=2,3$,
  that cannot be (formally) lifted to characteristic $0$.
\end{theorem}

  \begin{quote}
    {\bf Question:}
    Are there non-liftable Calabi-Yau 3-fold
    in $p\geq 5$?
  \end{quote}

  By the time when the author writes this article,
we only know  non-liftable 3-dimensional Calabi-Yau spaces,
(i.e., algebraic spaces which are not schemes)
in the case $p\geq 5$.
Moreover, whether there exist non-liftable CY $n$-folds for $n\geq 4$
  is unclear.

\subsection{Hodge decomposition and  Kodaira-Akizuki-Nakano vanishing}

Let 
\begin{equation*}
 F: X\To X^{(p)}
  \end{equation*}
be a relative Frobenius morphism.
For a complex $L^\bullet$ and $n\in\ZZ$ we denote by $T^\bullet:=\tau_{<n}L^\bullet$
the complex such that
\begin{equation*}
  T^i =
  \left\{
  \begin{array}{ll}
    L^i & \mbox{for }i\leq n-2\\
    \Ker (d: L^{n-1}\To L^n) & \mbox{for } i = n-1\\
    0  & \mbox{for }i\geq n
    \end{array}
  \right.
  \end{equation*}
Also, we denote by $W_2$ the ring of second Witt vectors
$W(k)/p^2W(k)$, where $W(k)$ is the ring of Witt vectors over $k$.

\begin{theorem}[Deligne-Illusie \cite{DI87}]
  \label{thm:DI87}
  Let $k$ be a perfect field of $\chara(k)=p>0$ and $X$ a smooth scheme
  over $k$. If $X$ is liftable over $W_2$, then we have
  \begin{equation*}
\varphi: \Dirsum_{i<p}\Omega^i_{X^{(p)}/k}[-i] \overset{\iso}{\To}
    \tau_{<p}F_*\Omega_{X/S}^\bullet
  \end{equation*}
  which is an isomorphism in the derived category $D(X^{(p)})$
  of $\OO_{X^{(p)}}$-modules with $F$ action
  such that
  \begin{equation*}
    {\cal H}^i\varphi = C^{-1}:
    \Dirsum_i\Omega^i_{X^{(p)}/k}\To \Dirsum_i {\cal H}^iF_*\Omega^\bullet_{X/k}
  \end{equation*}
  for $i<p$, where $C$ is the  Cartier operator.
\end{theorem}

From this theorem, we obtain the following consequences.

\begin{corollary}[Hodge decomposition]
  Let $X$ be a smooth proper scheme over a perfect field $k$ of $\chara(k)=p>0$
  and $\dim \leq p$.
  If $X$ is liftable over $W_2$, 
  then Hodge to de Rham spectral sequence
  \begin{equation*}
    E_1^{pq}=H^q(X,\Omega^p) \Rightarrow H^{p+q}_{DR}(X/k)
  \end{equation*}
  degenerates at $E_1$.
  In particular, we have Hodge decomposition:
  \begin{equation*}
    H_{DR}^n(X/k)\iso \Dirsum_{i=0}^n H^{n-i}(X,\Omega^i)
  \end{equation*}
  for $n\in\ZZ$.
\end{corollary}

\begin{corollary}[Kodaira-Akizuki-Nakano vanishing]
  \label{KNAposchar}
    Let $X$ be a smooth projective scheme over a perfect field $k$ of
    $\chara(k)=p>0$ and $L$ an ample line bundle on $X$.
    If $X$ is liftable over $W_2$, then we have
    \begin{equation*}
      H^j(X, \Omega^i\tensor L^{-1})=0
      \qquad    \mbox{for } i+j<\inf(\dim X, p).
    \end{equation*}
    In particular, if $\dim X\leq p$, we have
    Kodaira vanishing $H^i(X, L^{-1})=0$ for $i<\dim X$.
  \end{corollary}

For Calabi-Yau 3-folds, the following question is widely open,
apart from partial answers.

\begin{quote}
{\bf Question:} For a non-liftable Calabi-Yau 3-fold $X$, is it
liftable over $W_2$? If not, does Kodaira vanishing hold?
\end{quote}

\subsection{Supersingularity}

For a Calabi-Yau $n$-fold $X$ ($n\geq 2$), Artin-Mazur functor
\cite{AM}
\begin{equation*}
  \Phi_X^n: {\rm Art}\To {\rm Abgr}
\end{equation*}
from the category of Artinian
local rings ${\rm Art}$ to the category of abelian groups
${\rm Abgr}$ is defined by
\begin{equation*}
  \Phi_X^n(S) := \Ker (H_{et}^n(X\times_k S, {\Bbb G}_m)\To
  H^n_{et}(X, {\Bbb G}_m))
\end{equation*}
and this is pro-representable by a 1-dimensional formal group
scheme $M$. Namely, we have 
\begin{equation*}
  \Phi_X^n(-) = \Hom(-,M).
\end{equation*}
It is known that a 1-dimensional formal group scheme $M$
is the formal additive group scheme
$\hat{{\Bbb G}}_a$ or a $p$-divisible formal group scheme.
The group operation of a formal group scheme can be described
by a formal group law $F(X,Y)\in k[[X,Y]]$ and we can define
the {\em height} of the formal group law. By definition, {\em height} $h(X)$
of a Calabi-Yau variety $X$ is the height of the formal group law.
See \cite{Haz} for the detail of formal group law for the group scheme $M$.

The height $ht(X)$ has convenient description in terms of Serre
cohomologies $H^*(X, W\OO_X)$ \cite{Serre} and crystalline
cohomologies $H_{cris}^*(X/W)$ \cite{Ber,BO}, from which we deduce
characterization of supersingular Calabi-Yau varieties.

We will denote by $K$ the quotient field of $W:=W(k)$.

\begin{proposition}
  For a Calabi-Yau variety $X$ of dimension $n=\dim X$ over an algebraically
  closed field $k$,  we have 
  \begin{eqnarray*}
    \lefteqn{ht(X)}\\
    &=& \left\{
    \begin{array}{ll}
      \dim_K H^n(X, W\OO_X)\tensor_W{K}
      & \\
      \quad = \dim_K (H^n_{cris}(X/W)\tensor_W{K})_{[0,1)}\; (<\infty) &
        \mbox{if }H^n(X,W\OO_X)\tensor_W{K}\ne{0}\\
        & \\
      \infty & \mbox{if }H^n(X,W\OO_X)\tensor_W{K}={0}\\
      \end{array}
    \right.
  \end{eqnarray*}
  where $(-)_{[0,1)}$ denotes the $K$-vector subspace of slopes
    between $0$ and $1$.
\end{proposition}

\begin{definition}[ordinary/supersingular]
  Let $X$ be a Calabi-Yau variety of $\dim X=n(\geq 2)$.
  Then, we call $X$ is {\em ordinary} if $h(X)<\infty$ and
  {\em supersingular} if $h(X)=\infty$.
  \end{definition}

\begin{proposition}[cf.\cite{A74}]
  If $n=2$, i.e., $X$ is a K3 surface,
  we have  $1\leq h(X)\leq 10$ or $h(X)=\infty$.
  The supersingular case $h(X)=\infty$ exists only in positive characteristic.
\end{proposition}

On only K3 surfaces but also for
higher dimensional Calabi-Yau varieties, supersingularity has
close relation with unique feature of geometry in positive
characteristic.

For a Calabi-Yau variety $X$ in characteristic $0$,
the Betti number $b_n(X)$, $n=\dim X$, can never be trivial since
by Hodge decomposition $b_n(X)\geq \dim_k H^0(X, \Omega_X^n)
= \dim_k H^0(X, \OO_X)=1$. In positive characteristic, we consider 
the \'etale Betti numbers 
\begin{equation*}
b_i(X):=\dim_{\QQ_\ell} H^i_{et}(X, \QQ_\ell)
:= \varprojlim H^i_{et}(X, \ZZ/\ell^r\ZZ)\tensor_{\ZZ_\ell}\QQ_\ell
\end{equation*}
for a prime number $\ell(\ne p)$ (see for example \cite{Milne}).

\begin{proposition}[cf.Prop.~8.1~\cite{Schr04} ]
  \label{bettizeroSS}
Calabi-Yau $n$-fold with $b_n(X)=0$ is supersingular.
\end{proposition}

Calabi-Yau 3-folds with trivial 3rd Betti number
do exist as will be presented
in section~\ref{sec:3rdlec}.
We note that the converse of Prop.~\ref{bettizeroSS}
does not hold since we have $b_2(X)=22$
for a supersingular K3 surface $X$.\\

Recall that a variety  $X$ is called {\em uniruled} (or {\em unirational}) if
there exists a dominant rational morphism $\varphi : Y\times \PP^1 \to X$
(or $\varphi: \PP^n\to X$) for some variety $Y$.
In characteristic~$0$, a Calabi-Yau variety cannot be uniruled since
otherwise $\omega_X$ cannot be tirivial because of the following fact:

\begin{theorem}[\cite{MM86}]
  Let $X$ be a smooth projective variety over $\CC$. Then
  $X$ is unuruled if and only if
  there exists a non-empty open subset $U\subset X$
  such that for all $x\in U$ there exists an irreducible
  curve $C$ through $x$ with $(K_X,C)<0$.
\end{theorem}

However, uniruled Calabi-Yau 3-folds do exist as will be presented
in section~\ref{sec:3rdlec}.

\begin{proposition}[Theorem~1.3 \cite{H99}, Theorem~3.1~\cite{H01}]
  \label{uniruledsupersingular}
  A uniruled Calabi-Yau $n$-fold $X$ is supersingular.
\end{proposition}

The converse to Prop.~\ref{uniruledsupersingular} is
open for $n\geq 3$. For $n=2$, C.~Liedtke  showed that
a supersingular K3 surface is unirational \cite{Lied15}.

\begin{quote}
  {\bf Question:} Is there any relation between non-liftability
  to characteristic $0$ and supersingularity?
\end{quote}

\section{Deformation theory of CY 3-folds in positive characteristic}

\subsection{Infinitesimal Lifting}

Let $X$ be a smooth projective variety over a perfect field $k$ of
$\chara(k)=p>0$
and $(A,\mm)$ a complete Noetherian local domain of mixed characteristic,
i.e., $A = \varprojlim A/\mm^n$, $\chara(A)=0$ and $k=A/\mm$.
We denote the quotient field of $A$ by $K$.

\begin{example}
  A typical situation we have mostly in mind is that $A$ is
  the ring of Witt vectors $(W(k), pW(k))$ over $k$
  or finite extension of $W(k)$.
\end{example}

\subsubsection{formal spectrum $\operatorname{Spf}(A)$ }
Set $A_n:= A/\mm^{n+1}$ for $n\geq 0$.
They are Artinian local rings and in particular $A_0= k$.
Also set $S_n:= \Spec A_n$.
Then we have an increasing sequence of infinitesimal
neighborhoods $S_n\subset S_{n+1}$, but all the $S_n$ have the same 
underlying space, which we denote by
$\operatorname{Spf}(A)$. We define its structure sheaf as 
\begin{equation*}
  \OO_{\operatorname{Spf(A)}}:= \varprojlim \OO_{S_n}
\end{equation*}
We have $\Gamma(\operatorname{Spf}(A), \OO_{\mathcal S}) = A$.

\subsubsection{formal lifting via infinitesimal lifting}

A lift ${\cal X}$ of $X$ over $A$ is a scheme flat over $A$,
${\mathcal X}\to \Spec(A)$, such that $X \iso {\mathcal X}\times_A k$.
Its closed fiber (or special fiber)
$X$ in characteristic $\chara(k)=p>0$ is lifted to
the generic fiber ${\cal X}_\eta := {\mathcal X}\times_AK$
in $\chara(K)=0$.

\begin{equation*}
    \begin{array}{ccccc}
X          &\hookrightarrow & {\cal X}   & \hookleftarrow & {\cal X}_\eta \\
\downarrow &                & \downarrow &                & \downarrow  \\
\Spec(k)   &\hookrightarrow & \Spec(A)   & \hookleftarrow & \Spec(K)\\
\end{array}
\end{equation*}

One way to construct such a lifting ${\cal X}$ is {\em infinitesimal lifting}
of $X$. Consider the  short exact sequence:
\begin{equation}
  \label{smallext}
  0\To \mm^n/\mm^{n+1} \To A_n \To A_{n-1}\To 0
\end{equation}
We note that $\mm^n/\mm^{n+1} \iso \mm \tensor k$ is a finite $k$-vector space.
Infinitesimal lifting is to try to lift $X$ over $A_0=k$ to a scheme
$X_1$ flat over $A_1=A/\mm^2$, then lift $X_1$ to a scheme $X_2$ flat over $A_2$
and so on to obtain a family $\{X_n\}_{n\geq 0}$ of schemes
such that $X_n\tensor_{A_n}A_{n-1} = X_{n-1}$ ($n\geq 1$).
Such a family $\{X_n\}$ is called a {\em formal family of deformations of $X$
  over $A$}. From this family, we obtain a formal lifting:

\begin{proposition}
  \label{prop21.1:HDT}
  Given such a family $\{X_n\}_{n\geq 0}$,
  we obtain  a Noetherian formal scheme $\tilde{X}$ flat over
  $\operatorname{Spf}(A)$ such that $X_n\iso \tilde{X}\times_A A_n$
  for all $0\leq n\in\ZZ$.
  \end{proposition}
\begin{proof}
  Define $\tilde{X}$ to be the locally ringed space formed by taking
  the topological space $X_0$ together with the sheaf of rings
  $\OO_{\tilde{X}} := \varprojlim \OO_{X_n}$. See, for example,
  Prop.~21.1~\cite{HDT}
  for the rest of proof.
\end{proof}

The formal scheme $\tilde{X}$ as in Prop.~\ref{prop21.1:HDT}
is called a {\em formal lifting} of $X$. A formal lifting
scheme is locally isomorphic to the completion of the closed fiber,
but it is not always so globally.
On the other hand, we say that $X$ can be {\em projectively lifted}
over the field $K$ of characteristic $0$ if there exists a projective
scheme ${\cal X}$ flat over $A$ with $X\iso {\cal X}\tensor_A k$.\\

Then we have questions: given a projective variety $X$ over $k$,
\begin{description}
\item [Q1:] when do we have a formal lifting $\tilde{X}$?
\item [Q2:] given a formal lifting $\tilde{X}$,
  when do we have a projective lifting
  ${\mathcal X}$, whose completion along the closed fiber
  is $\tilde{X}$?
\end{description}

\subsubsection{obstruction to formal lifting}

\begin{definition}[torsor/principal homogeneous-space]
  Suppose that a group $G$ acts on a non-empty set $S$.
    Then $S$ is called a {\em torsor} or a {\em principal homogeneous space}
    under the action of $G$
    if there exists one (and hence all) element $s_0\in S$ such that
    $G\iso S$ via $g\mapsto g(s_0)$
\end{definition}

Now the answer to the question {\bf Q1} is that
if $H^2(X,T_X)=0$ we always have a formal lifting of $X$:

\begin{proposition}
\label{cor10.3:HDT}
  The obstruction to lifting at each step lie in $H^2(X,T_X)$.
  If a lifting exists, the set of equivalence classes of
  all such is a torsor
  under $H^1(X,T_X)$, which means that
  there is a one-to-one correspondence between the set of
  equivalence classes of liftings
  and $H^1(X,T_X)$.
\end{proposition}
\begin{proof}
  Consider the short exact sequence (\ref{smallext}) and assume that
  we have a scheme $X_{n-1}$ over $A_{n-1}$.  Then, we can show that
  there is just one obstruction in $H^2(X, T_X\tensor
  \mm^n/\mm^{n+1})$ for the existence of a lifting $X_n$ of $X_{n-1}$
  over $A_n$.  Since $\mm^n/\mm^{n+1}$ is a finite dimensional vector
  space over $k$, we know that $H^2(X,T_X)$ implies no obstruction for
  the lifting.  See for example Cor.10.3~\cite{HDT} for the rest of
  the proof.
  \end{proof}

Here we define rigidity of scheme, which will be used later.

\begin{definition}[rigid scheme]
  A scheme $X$ is called {\em rigid} if all of its deformations over
  the dual numbers are trivial, i.e., $X\times_kD$, where
  $D:=\Spec k[\epsilon]/(\epsilon^2)$.
  \end{definition}

\begin{proposition}
  For a smooth scheme $X$ over $k$,
  the following are equivalent:
  \begin{enumerate}
  \item $X$ is rigid;
  \item all the deformations over an Artinian $k$-algebra $A$
    are trivial, i.e., $X\times_k\Spec A$
  \item $H^1(X,T_X)=0$.
    \end{enumerate}
  \end{proposition}
\begin{proof}
  See Excercise~10.4 and Theorem~5.3.
  \end{proof}

\subsubsection{obstructions to projective lifting}

The answer to the question {\bf Q2} is

\begin{theorem}
  \label{HDT:th22.1}
  In the situation described above, 
  assume that $H^2(X,\OO_X)=0$ and $H^2(X, T_X)=0$.
  Then, $X$ can be projectively lifted to
  characteristic $0$.
\end{theorem}
\begin{proof}
  Since $H^2(X,T_X)=0$, we have an infinitesimal lifting $\tilde{X}$
  by Prop.~\ref{cor10.3:HDT}. The obstruction to lifting
  an invertible sheaf lies in $H^2(X,\OO_X)$ (see, for example,
  Theorem~6.4(a)~\cite{HDT}), which is trivial by
  assumption. Thus an ample line bundle $L$ of $X$ can be
  lifted to $\tilde{X}$. Then applying Grothendieck's
  algebraization theorem we obtain a projective variety ${\mathcal X}$
  whose completion along
  the closed fiber $X$ is the formal lifting $\tilde{X}$.
  See  Theorem.~22.1~\cite{HDT} for the detail.
\end{proof}

\begin{corollary}
For a Calabi-Yau $n$-fold $(n\geq 3)$,
existence of a formal lifting 
implies existence of a projective lifting.
\end{corollary}

Thus, projective lifting problem is a little simpler for $\dim X\geq 3$.
Compare with the comments after Theorem~\ref{DI:K3}.

\subsection{$W_2$-liftability of ordinary CY $n$-folds}

Although there exist supersingular Calabi-Yau 3-folds
which cannot be lifted over $W(k)$, or even over $W_2$,
the situation is a little different for ordinary Calabi-Yau varieties.
Namely,  F.~Yobuko \cite{Yobu} proved recently that ordinary CY $n$-folds
always lift over $W_2$, so that in particular, for $\chara(k)\geq 3$,
Kodaira vanishing holds by Cor.~\ref{KNAposchar}. Precisely,

\begin{theorem}[Yobuko~\cite{Yobu}]
  \label{thm:yobuko}
  Let $X$ be a Calabi-Yau variety over an algebraically closed field
  $k$ of $\chara(k)=p>0$ with $ht(X)<\infty$. Then $X$ admits a flat
  lift over $W_2$.
  \end{theorem}

Hence by Theorem~\ref{thm:DI87}  we immediately obtain

\begin{corollary}
  Let $X$ be an ordinary Calabi-Yau variety over an algebraically closed
  field $k$ with $\dim{X}\leq p=\chara(k)$. Then
  Kodaira vanishing holds for $X$.
\end{corollary}

The key ideas of the proof are quasi-Frobenius splitting and
splitting height. To understand these notions, we need (iterated)
Cartier operator, which we will explain now.

\subsubsection{Cartier operator}
Let $X$ be a scheme with $p\OO_X=0$ with the structure morphism
$\pi: X\to S$ as ${\Bbb F}_p$-schemes. Then we consider
the absolute Frobenius morphisms $F_X$ and $F_S$ together with 
the relative Frobenius morphism $F_{X/S}: X\To X^{(p)}$:
\begin{equation*}
  \begin{CD}
    X @<{W}<<   X^{(p)} @<{F_{X/S}}<<   X \\
    @V{\pi}VV    @V{\pi^{(p)}}VV      @V{\pi}VV \\
    S @<{F_S}<< S  @=                 S\\
  \end{CD}
  \qquad
  \begin{array}{l}
    X^{(p)}:= X\times_S (S,F_S)\\
    F_X = W\circ F_{X/S} 
    \end{array}
\end{equation*}
Notice that since $\Omega_{X^{(p)}/S}^i = W^*\Omega_{X/S}^i$,
we have $W^*(df)\in \Omega_{X^{(p)}/S}^1$ for $df\in \Omega_{X/S}^1$.\\

\begin{example}
  If $S = \Spec R$ and $X = \Spec R[\overline{T}]$
  , $R[\overline{T}]:=R[T_1,\ldots, T_n]$,
then 
$W: X^{(p)}\to X$ is induced by the ring homomorphism
\begin{equation*}
  W^\sharp: R[T_1,\ldots, T_n]
  \To R[T_1,\ldots, T_n]\tensor_R(R, F_R)
  \end{equation*}
such that $W^\sharp(T_i) = T_i$ for all $i$ 
and $W^\sharp(f)= F_R(f)=f^p$  for all $f\in R$.
Also, $F_{X/S}$ is induced from the ring homomorphism
\begin{equation*}
  F_{X/S}^\sharp: R[T_1,\ldots, T_n]\tensor_R(R, F_R)
  \To R[T_1,\ldots, T_n]
\end{equation*}
such that $F_{X/S}^\sharp(T_i) = T_i^p$ and $F_{X/S}^\sharp(f) =f$
for $f\in R$.
Then
\begin{eqnarray*}
  \lefteqn{W^*\Omega^i_{X/S}}\\
    &= &(R[\overline{T}], F_R)\tensor_{R[\overline{T}]} \Omega_{X/S}^i\\
  &= &(R[\overline{T}], F_R)\tensor_{R[\overline{T}]}
\left\{\sum_{1\leq j_1<\cdots, j_i\leq n}
r_{j_1\cdots j_i}\cdot dT_{j_1}\wedge\cdots\wedge dT_{j_i}
\;\vert\; r_{j_1\cdots j_i}\in R[\overline{T}]\right\}\\
 & = & \Omega^i_{X^{(p)}/S}
\end{eqnarray*}
\end{example}

Next 
we consider the de Rham complex $(\Omega_{X/S}^\bullet, d)$
and let 
  \begin{equation*}
    Z_1\Omega_{X/S}^i:= \Ker (d: \Omega_{X/S}^i\To \Omega_{X/S}^{i+1})
    \quad\mbox{and}\quad
    B_1\Omega_{X/S}^i:= \Im (d: \Omega_{X/S}^{i-1}\To \Omega_{X/S}^{i}).
    \end{equation*}
Since 
 $d: \Omega_{X/S}^i\to \Omega_{X/S}^{i+1}$ is a $\OO_{X^{(p)}}$-linear map
for every $i\geq 0$, $Z_1\Omega_{X/S}^i$, $B_1\Omega_{X/S}^i$
  and the cohomology sheaf
  ${\mathcal H}^*(\Omega_{X/S}^\bullet)$ are all $\OO_{X^{(p)}}$-modules.
  We sometimes write as ${\mathcal H}^i((F_{X/S})_*\Omega_{X/S}^\bullet)$ to
  stress this fact.

Now we define the Cartier operator.

\begin{theorem}
  \label{thm:cartieroperator}
  There exists a morphism
  \begin{equation*}
    C^{-1}=C_{X/S}^{-1}: \Omega_{X^{(p)}/S}^i \To {\cal H}((F_{X/S})_*\Omega_{X/S}^\bullet)
  \end{equation*}
  such that
  \begin{enumerate}
  \item $C^{-1}(1)=1$,
  \item $C^{-1}(\omega\wedge\tau) = C^{-1}(\omega)\wedge C^{-1}(\tau)$
  \item $C^{-1}(W^*(df)) = \mbox{the cohomology class of }f^{p-1}df$
  for $f\in \OO_X$.
  \end{enumerate}
  $C^{-1}$ is isomorphic if $X/S$ is smooth.
\end{theorem}
\begin{proof}
  See, for example, Theorem.~7.1~\cite{Katz70}.
  \end{proof}
    
\begin{definition}[Cartier operator]
  Assume that $X/S$ is smooth. Then the morphism $C_{X/S}:Z_1\Omega^\bullet_{X/S}
  \To \Omega^\bullet_{X^{(p)}/S}$ induced by Theorem~\ref{thm:cartieroperator}
  \begin{equation}
    \label{cartierseq}
    0\To B_1\Omega_{X/S}^\bullet \To Z_1\Omega^\bullet_{X/S}
    \overset{C_{X/S}}{\To} \Omega_{X^{(p)}/S}^\bullet
    \To 0
    \end{equation}
  is called the {\em Cartier operator}.
  \end{definition}

\subsubsection{iterated Cartier operator and quasi Frobenius splitting}
We recall the iterated Cartier operator defined in \cite{Ill79}
(see also \cite{KvdG}).
Let $X$ be a smooth scheme over a scheme $S$ and $n\in{\Bbb N}$. Then we
define
\begin{eqnarray*}
  X^{(p^n)} &:=&  X \times_S(S, F_S^n) = (X^{(p^{n-1})})^{(p)} \\
F_{X/S}^n &:= & F_{X^{(p^{n-1})}/S}\circ
F_{X^{(p^{n-2})}/S}\circ \cdots \circ F_{X/S}
=  (F_{X/S}^{n-1})^{(p)}\circ F_{X/S}
\end{eqnarray*}

\begin{equation*}
  \begin{CD}
    X @<{W^n}<<   X^{(p^n)} @<{F^n_{X/S}}<<   X \\
    @V{\pi}VV    @V{\pi^{(p^n)}}VV      @V{\pi}VV \\
    S @<{F^n_S}<< S  @=                 S\\
  \end{CD}
  \qquad
    F^n_X = W^n\circ F^n_{X/S} 
\end{equation*}\\

where the right square is expanded as follows:

\begin{equation*}
  \begin{CD}
    X^{(p^n)} @<{F_{X^{(p^{n-1})}/S}}<<   X^{(p^{n-1})} @<{F_{X^{(p^{n-2})}/S}}<<   X^{(p^{n-2})}
            @<<<\cdots @<<<    X^{(p)}   @<{F_{X/S}}<< X\\
            @V{\pi^{(p^n)}}VV          @V{\pi^{(p^{n-1})}}VV       @V{\pi^{(p^{n-2})}}VV
            @.  @V{\pi^{(p)}}VV   @V{\pi}VV \\
    S  @=                 S  @=  S @= \cdots @= S @= S\\
  \end{CD}
\end{equation*}

Then we define
\begin{equation*}
  \begin{array}{ccc}
    0\subset B_1\Omega_{X/S}^i \subset \cdots
    &\subset B_n\Omega_{X/S}^i \subset
   B_{n+1}\Omega_{X/S}^i \subset  \cdots & \\
   & \cdots \subset Z_{n+1}\Omega_{X/S}^i 
   \subset Z_n\Omega_{X/S}^i \subset &\cdots \subset
  Z_1\Omega_{X/S}^i \subset \Omega_{X/S}^i
  \end{array}
\end{equation*}
as follows:
\begin{eqnarray*}
 B_0\Omega^i_{X/S}&=& 0,\quad Z_0\Omega^i_{X/S} = \Omega_{X/S}^i\\
  B_n\Omega_{X^{(p)}/S}^i & \overset{C^{-1}_{X/S}}{\iso}& B_{n+1}\Omega_{X/S}^i/B_1\Omega_{X/S}^i\\
  Z_n\Omega_{X^{(p)}/S}^i & \overset{C^{-1}_{X/S}}{\iso}& Z_{n+1}\Omega_{X/S}^i/B_1\Omega_{X/S}^i\\
\end{eqnarray*}
Then we have

\begin{proposition}[(0.2.2.5) and (0.2.2.6.3) \cite{Ill79}]
  \begin{equation*}
     \Omega_{X^{(p^n)}/S}^i
     \overset{C_{X/S}^{-n}}{\iso} Z_n\Omega_{X/S}^i/B_n\Omega_{X/S}^i
     \qquad\mbox{and}\quad
    {\cal H}^i(\Omega_{X^{(p^n)}/S}^\bullet)
     \overset{C_{X/S}^{-n}}{\iso}  Z_{n+1}\Omega_{X/S}^i/B_{n+1}\Omega_{X/S}^i
  \end{equation*}
where, we write $C_{X^{(p^i)}/S}$, $i=0,1,2,\ldots$,  simply as $C_{X/S}$.
\end{proposition}
\begin{proof}
  We have, for  $m=n$ or $n+1$,
  \begin{eqnarray*}
    Z_m\Omega_{X/S}^i/B_m\Omega_{X/S}^i
    &\iso &
    \frac{Z_m\Omega_{X/S}^i/B_1\Omega_{X/S}^i}
    {B_m\Omega_{X/S}^i/B_1\Omega_{X/S}^i}\\
    &\overset{C_{X/S}}{\iso}&
    Z_{m-1}\Omega_{X^{(p)}/S}^i/B_{m-1}\Omega_{X^{(p)}/S}^i
    \iso
        \frac{Z_{m-1}\Omega_{X^{(p)}/S}^i/B_1\Omega_{X^{(p)}/S}^i}
    {B_{m-1}\Omega_{X^{(p)}/S}^i/B_1\Omega_{X^{(p)}/S}^i}\\
    &\overset{C_{X/S}}{\iso}&\cdots\\
    &\overset{C_{X/S}}{\iso}&
    Z_{1}\Omega_{X^{(p^{m-1})/S}}^i/B_{1}\Omega_{X^{(p^{m-1})}/S}^i
    = {\cal H}^i(\Omega_{X^{(p^{m-1})}/S}^\bullet)\\
    &\overset{C_{X/S}}{\iso}&
    \Omega_{X^{(p^{m})}/S}^i.
  \end{eqnarray*}
  \end{proof}

\begin{remark}
  Notice that we can also define
    $B_{n+1}\Omega_{X/S}^i = C_{X/S}^{-1}(B_{n}\Omega_{X/S}^i)$
    and
    $Z_{n+1}\Omega_{X/S}^i = \Ker (dC_{X/S}^{n})$
    as in  \cite{KvdG}.
\end{remark}

From now on, we set $S = \Spec(k)$ and we only consider the absolute Frobenius
morphism.
First of all we consider the exact sequence
\begin{equation*}
  0\To \OO_X\To (F_X)_*\OO_X \overset{d}{\To} B_1\Omega_{X}^1\To 0\qquad (e)_1
\end{equation*}

\begin{definition}[Frobenius splitting \cite{MR85}]
  $X$ is said to be {\em Frobenius split} if $(e)_1$ splits.
\end{definition}

We will now extend this notion: By pulling back along (the surjection)
$B_m\Omega_X^1\overset{C^{m-1}}{\To} B_1\Omega_X^1$, we obtain the
following diagram:
\begin{equation}
  \label{emtoe1}
  \begin{CD}
  0 @>>>  \OO_X @>>>(F_X)_*\OO_X @>{d}>> B_1\Omega_{X}^1 @>>> 0\qquad (e)_1\\
  @.      @|         @AAA                @AA{C^{m-1}}A        @. \\
  0 @>>>  \OO_X @>>> {\cal F}_m  @>>> B_m\Omega_{X}^1 @>>> 0\qquad (e)_m\\
  \end{CD}
\end{equation}
where  ${\cal F}_m$ is the extension of $B_m\Omega_{X}^1$ by
$\OO_X$ fitting in this diagram.

Now we define
\begin{definition}[splitting height and quasi Frobenius split]
  We define the {\em splitting height} $\operatorname{sht}(X)$ of $X$ by
  \begin{equation*}
     \operatorname{sht}(X) := \min\{m\;\vert\; (e)_m\mbox{ splits}\}
  \end{equation*}
  $X$ is called {\em quasi Frobenius split} if
  $\operatorname{sht}(X)<\infty$.
  \end{definition}

The splitting height gives a new interpretation of height for Calabi-Yau varieties.

\begin{proposition}[Prop.~6~\cite{Yobu}]
  \label{prop6:yobu}
For a Calabi-Yau variety $X$ of $\dim{X}=n$, we have $ht(X)=sht(X)$
\end{proposition}
\begin{proof}
  By the following Lemma~\ref{ht1CY}, $ht(X)=1$ if and only if
  $sht(X)=1$. Thus we may regard $ht(X)\geq 2$ and
  $sht(X)\geq 2$ and $(e)_1\ne 0$.
  From the exact sequence
  \begin{equation*}
    0\To B_{m-1}\Omega_X^1\To
    B_m\Omega_X^1\overset{C^{m-1}}{\To} B_1\Omega_X^1\To 0
  \end{equation*}
  we have a long exact sequence
  \begin{equation*}
    \Ext^1(B_1\Omega_X^1, \OO_X)
    \overset{(C^{m-1})^*}{\To}
    \Ext^1(B_m\Omega_X^1, \OO_X)
    \overset{\varphi_m}{\To}
    \Ext^1(B_{m-1}\Omega_X^1, \OO_X)
    \To \Ext^2(B_1\Omega_X^1,\OO_X).
  \end{equation*}
  By the definition of $(e)_m$, we have $(C^{m-1})^*((e)_1) = (e)_m$.

  Now recall Prop.~3.1\cite{KvdG}, i.e.,
  \begin{equation*}
    \dim H^i(X, B_m\Omega_X^1)
    = \left\{
    \begin{array}{ll}
      0   &  (i\ne n-1, n)\\
      \min\{m, ht(X)-1\} & (i=n-1, n).
      \end{array}
    \right.
  \end{equation*}
  and the fact that
  $B_m\Omega_X^1$ is locally free as a $\OO_{X^{(p^m)}}$-module
  (see for example. Prop.~0.2.2.8(a)~\cite{Ill79}). Then by Serre duality,
  we compute
   $\Ext^2(B_1\Omega_X^1,\OO_X)
  = \Ext^2(\OO_X, (B_1\Omega_X^1)^\vee) =
  H^2(X, (B_1\Omega_X^1)^\vee)
  = H^{n-2}(X, B_1\Omega_X^1)^\vee =0$,
  which implies that $\varphi_m$ is surjective.
  Similarly
  $\dim \Ext^1(B_m\Omega_X^1,\OO_X) = \min\{m, ht(X)-1\}$.
  Thus we know in particular that $\Ext^1(B_1\Omega_X^1,\OO_X)$
  is $1$-dimensional and generated by $(e)_1$,
  which means that $\Ker\varphi_m = \Im (C^{m-1})^*$ is generated by
  $(e)_m$. Consequently, we know that
  $(e)_m=0$ if and only if
  $\varphi_m$ is isomorphic.
  
  Now let $s:=sht(X)$. Then $(e)_m\ne{0}$ for all $m<s$. Thus
  if we set $d_m:= \dim \Ext^1(B_m\Omega_X^1,\OO_X)
  = \min\{m, ht(X)-1\}$, we have that $\varphi_m$
  is non-isomorphic for $m<s$ and $\varphi_s$ is isomorphic.
  Hence we have
  \begin{equation*}
    d_s=d_{s-1}>d_{s-2}>\cdots > d_2 > d_1=1
  \end{equation*}
  and then
  $d_s=\min\{s, ht(X)-1\}=\min\{s-1, ht(X)-1\}=d_{s-1}\geq s-1$,
  from which we obtain $ht(X)= sht(X)$.
  \end{proof}

The following result is well-known to experts.

\begin{lemma}
  \label{ht1CY}
  For a Calabi-Yau $n$-fold $X$, $ht(X)=1$ if and only if
  $X$ is Frobenius-split.
\end{lemma}
\begin{proof}
  By \cite{KvdG}, we have
  \begin{equation*}
    ht(X) = \min\{i\geq 1\;\vert\; F: H^n(X, W_i\OO_X)\to H^n(X, W_i\OO_X)
    \mbox{ is non-trivial}\}
  \end{equation*}
  so that $ht(X)=1$ is equivalent
  to $F:H^n(X,\OO_X)\to H^n(X, \OO_X)$ being non-trivial.

  On the other hand,
   $X$ being Frobenius-split means that
  \begin{equation*}
    0\To \OO_X \To F_*\OO_X \overset{d}{\To} B_1\Omega_X^1\To 0
  \end{equation*}
  splits. Taking the dual $(-)^\vee = \Hom(-,\OO_X)$
  it is equivalent that
  \begin{equation*}
    0\To (B_1\Omega^1_X)^\vee \To (F_*\OO_X)^\vee \To
    \Hom(\OO_X,\OO_X)(\iso \OO_X)\To 0
  \end{equation*}
  splits, where we note that, since $X$ is smooth, all these sheaves
  are locally free.
  This means that the identy $id: \OO_X\to \OO_X$ can be
  lifted to $F_*\OO_X\to \OO_X$, so that it is equivalent that
  \begin{equation*}
    H^0(X, (F_*\OO_X)^\vee) \To H^0(X, \OO_X)
  \end{equation*}
  is non-zero. By Serre-duality and $\omega_X\iso \OO_X$, this is equivalent
  with non-triviality of the Frobenius
  \begin{equation*}
     H^n(X,\OO_X)\To H^n(X, F_*\OO_X).
  \end{equation*}
\end{proof}

\subsubsection{Mehta-Srinivas deformation theory}

V.B.~Mehta and V.~Srinivas \cite{MS87,Sri90} formulated a simplified version of
Deligne-Illusie theory \cite{DI87} of $W_2$-lifting.
In this subsetction, we briefly summarize their results.

For any smooth variety $X$, the exact sequence $(e)_1$ induces the connecting homomorphism
\begin{equation}
  \label{eq:delta}
\Ext^1(\Omega_X^1, B_1\Omega_X^1) \overset{\delta}{\To} \Ext^2(\Omega_X^1, \OO_X).
\end{equation}
Then there is a class
\begin{equation*}
  obs(X) \in\Ext^2(\Omega_X^1, \OO_X)
\end{equation*}
whose non-vanishing is the obstruction to  lifting $X$ over $W_2$.
Furthermore, there is a class
\begin{equation*}
  obs(X,F_X)\in \Ext^1(\Omega_X^1, B_1\Omega_X^1)
\end{equation*}
whose non-vanishing is the obstruction to
lifting the pair $(X,F_X)$ to $W_2$.
This corresponds to the exact sequence
(\ref{cartierseq}).
Now the connecting homomorphism behaves like a forgetful map:
\begin{equation*}
  \delta(obs(X,F_X)) = obs(X).
\end{equation*}
Then we have 
\begin{proposition}
  If a smooth variety $X$ is Frobenius-split, then $X$ admits a
  lift over $W_2$.
\end{proposition}
\begin{proof}
  If $X$ is Frobenius-split, then $(e)_1$ splits so that $\delta=0$.
  Then we have $obs(X) = \delta(obs(X,F_X))=0$.
  \end{proof}

\subsubsection{Proof of Theorem~\ref{thm:yobuko}}

In view of Proposition~\ref{prop6:yobu}, we have only to show

\begin{proposition}[Prop.~5~\cite{Yobu}]
  Every quasi-Frobenius split variety admits a flat lift
  to $W_2$.
\end{proposition}
\begin{proof}
  Let $sht(X)=m<\infty$.
  From the diagram (\ref{emtoe1}) of 
  the extensions $(e)_m$ and $(e)_1$,
  we have the following commutative diagram:
  \begin{equation*}
    \begin{CD}
      \Ext^1(\Omega_X^1, B_1\Omega_X^1) @>{\delta}>> \Ext^2(\Omega_X^1, \OO_X) \\
      @A{(C^{m-1})_*}AA                                   @| \\
      \Ext^1(\Omega_X^1, B_m\Omega_X^1) @>{\delta_m}>> \Ext^2(\Omega_X^1,\OO_X).\\
      \end{CD}
  \end{equation*}
    Since $(e)_m=0$, we have $\delta_m=0$. Recall that
  $obs(X,F)\in \Ext^1(\Omega_X^1, B_1\Omega_X^1)$ and if
  it is in the image of $(C^{m-1})_*$, the commutativity of the
  diagram shows that $obs(X)=\delta(obs(X,F))=0$ and $X$ lifts over
  $W_2(k)$.
  But from the commutative diagram
  \begin{equation*}
    \begin{CD}
      0 @>>> B_1\Omega_X^1 @>>> Z_1\Omega_X^1 @>{C}>> \Omega_X^1 @>>> 0 \\
      @.   @A{C^{m-1}}AA       @A{C^{m-1}}AA           @|                \\
      0 @>>> B_m\Omega_X^1 @>>> Z_m\Omega_X^1 @>{C^m}>> \Omega_X^1 @>>> 0 \\        
      \end{CD}
    \end{equation*}
   we immediately know
     $obs(X,F)\in \Im (C^{m-1})_*$
  \end{proof}

The following problem seems to be still open:

  \begin{quote}
  {\bf Question:} 
  Is an ordinary Calabi-Yau $n$-fold $X$ liftable over
  $W_n$ ($n\geq 3$) or even liftable to characteristic $0$
  ?
  \end{quote}

  \section{Construction of non-liftable CY 3-folds
  \label{sec:3rdlec}}

In this section, after presenting what is known 
on the causes for non-liftability of Calabi-Yau varieties.
After that we overview the construction techniques
of the known examples of non-liftable Calabi-Yau threefolds.

\subsection{what causes non-liftability?}

The following two cases in which a Calabi-Yau variety is non-liftable
have been known so far: one is trivial highest Betti number
and the other is small resolution of liftable singular rigid
Calabi-Yau variety. We will explain them in the sequel and we still
have the following open question:

  \begin{quote}
    {\bf Question:}
    Are these two cases equivalent?
  \end{quote}

\subsubsection{non-liftability by $b_n(X)=0$ (due to Hirokado)}

By lifting to characteristic~$0$, we mean that ${\mathcal X}\to
\Spec(R)$ is proper and \underline{smooth} in this section.

\begin{proposition}
  \'{E}tale Betti numbers are preserved
  in lifting to characteristic $0$.
\end{proposition}
  \begin{proof}
    Let ${\mathcal X}\to S$ be a smooth and proper lifting of $X$
    with $X_\eta$ the generic fiber.
    Then, 
    by the proper and smooth base change theorem of \'etale cohomology,
    we have 
    \begin{equation*}
      H^i(X, \ZZ/\ell^n\ZZ) \iso H^i(X_\eta, \ZZ/\ell\ZZ)
    \end{equation*}
    for a prime $\ell(\ne p)$.
    Then by taking the limit, we obtain $b_i(X) = b_i(X_\eta)$.
    \end{proof}

\begin{corollary}
  \label{b3is0}
  For a Calabi-Yau $n$-fold $X$ in positive characteristic,
$\beta_{\dim{X}}(X)=0$ implies that $X$ is non-liftable to
  characteristic~$0$.
\end{corollary}
\begin{proof}
  Let $X_\eta$ be the lifting of $X$ to characteristic $0$.
  Since $X_\eta$ is also Calabi-Yau, the Hodge decomposition
  and Serre duality together with $\omega_X\iso \OO_X$
  show that
  \begin{equation*}
    b_n(X_\eta) = \sum_{i=0}^n\dim_K H^{n-i}(X_\eta, \Omega_{X_\eta}^i)
  \geq \dim_K H^{n}(X_\eta, \OO_{X_\eta})\iso H^0(X_{\eta},\OO_{X_\eta})^\vee
  =1.
  \end{equation*}
  Thus $0=b_n(X) = b_n(X_\eta)\geq 1$, a contradiction.
  \end{proof}

Notice that a Calabi-Yau $n$-fold with $b_n(X)=0$ is supersingular
by Prop.~\ref{bettizeroSS}.

\subsubsection{non-liftability by  mod $p$ reduction (due to Cynk and van Straten)
\label{section:cynk}}

We consider the situation in the following diagram:
\begin{equation*}
      \begin{CD}
      Y         @.                       @.              \\
      @V{\pi}VV               @.                  @. \\
      X         @>>>  {\mathcal X}       @<<<   {\cal X}_\eta \\
      @VVV                    @VVV                @VVV          \\
      \Spec(k)  @>>>          S=\Spec(R) @<<<   \Spec(K)
      \end{CD}
\end{equation*}
    where
$X$ is a rigid singular Calabi-Yau 3-fold
      with nodes as singularity,
${\mathcal X}$ is a lifting of $X$ over
      a complete Noetherian local domain $(R,\mm)$
      of mixed characteristic,
${\cal X}_\eta$ is the generic fiber, which is
      a (smooth) Calabi-Yau 3-fold over $K:=Q(R)$,
$\pi:Y\To X$ is a desingularization of
      nodes in $X$ with small resolution, which is locally described as
      \begin{equation*}
        X = \Spec k[[x,y,z,w]]/(xy-zw)
      \end{equation*}
      and
      \begin{equation*}
        \pi: Y= Bl_I(X)\To X
        \qquad (I=(x,z)\mbox{ or } I=(y,z))
        \end{equation*}
      is the blow-up centered at the ideal $I$.
Notice that the local analytic desingularization does not necessarily
extend to global algebraic desingularization in the category of
schemes: sometimes we have desingularization only in the
category of algebraic spaces.\\

    Then we have 
\begin{proposition}[Theorem~3.1 and 4.3~\cite{CyvS}]
        $Y$ is a non-liftable Calabi-Yau 3-fold.
\end{proposition}
\begin{proof}
  Assume that $Y$ has a lifting ${\mathcal Y}$ over $R$
  and set   $Y_n := {\mathcal Y}\times_R R/\mm^{n+1}$ for $1\leq n\in\ZZ$.
  Since $X$ has at most rational singularity, we can show
  that there exists commutative diagrams:
  \begin{equation*}
    \begin{array}{ccc}
      Y           &   \hookrightarrow &  Y_n  \\
      \downarrow  &                   & \downarrow \\
      X           &  \hookrightarrow  & X'_n
    \end{array}
  \end{equation*}
  for some $X'_n$ (blowing-down of lifting),
  but since $X$ is rigid, we must have $X'_n=X_n:=
  {\mathcal X}\times_R R/\mm^{n+1}$.
  Since $X$ has nodes, we have a section $\Spec R/\mm^{n+1}\to X_n$
  passing through a node. Since ${\mathcal X} = \varprojlim X_n$,
  ${\mathcal X}$ also has such a section and this means that
  ${\mathcal X}$ is singular, a contradiction.
  \end{proof}

\subsection{construction (I): quotient by foliation}

The first example of non-liftable Calabi-Yau 3-fold
\cite{H99} has been
constructed using quotient by foliation, which is a well-known
method to produce purely inseparable cover.

\subsubsection{foliation in positive characteristic}

Most of interesting pathologies in positive characteristic
are caused by purely inseparable extension of function fields.
Foliation (in algebraic setting) is a tool to produce such extensions
\cite{RS76, Ru84, Ek87, Jac75}.

The set of derivations $Der(K)$ has the structure of $p$-Lie algebra
(or $p$-closed Lie algebra, restricted Lie algebra), which means the
Lie algebra closed under Lie bracket and $p$th power.

\begin{definition}[exponent and $p$-basis]
  Let $K/k$ be an extension of fields of characteristic $p>0$.
  \begin{enumerate}
    \item [$(i)$] We say that $a\in K$ is of {\em exponent} $e(\geq 0)$
if $a$ is purely inseparable over $k$ and its minimal polynomial
is in the form of $x^{p^e}-c$ for some $c\in k-\{0\}$.
$K/k$ is called purely inseparable of {\em exponent $\leq e$}
if every element $a\in K$ is of exponent $\leq e$.
\item [$(ii)$] a set $B=\{\rho_1,\ldots, \rho_m\}\subset
  K$ is a {\em $p$-basis} of the extension if $K = k(K^p)(B)$,
  equivalently
$\{\rho_1^{k_1}\cdots\rho_m^{k_m}\;\vert\; 0\leq k_i<p\}$ is the
$k(K^p)$ basis of degree $p^m$ extension $K/k(K^p)$.
  \end{enumerate}
\end{definition}

For an extension $K/k$ of fields in characteristic $p>0$, we have $D:
K/k(K^p) \to K/k(K^p)$ for a derivation $D: K/k \To K/k$
over $k$.
Namely, we can consider $D$ to be a derivation over $k(K^p)$ and
$K/k(K^p)$ is of exponent $\leq 1$
since the extension $K/k(K^p)$ is described by $p$-basis.
Thus when we consider relation between derivation and purely inseparable
extension, we have only to consider the action of a derivation on 
$p$-basis of extension of exponent $\leq 1$.

Now let $K$ be a field of $\chara(K)=p>0$.
According to Jacobson's Galois theory for purely inseparable extensions
(of exponent one) \cite{Jac75}, there is a one-to-one correspondence
between
\begin{eqnarray*}
  {\cal E} & := &
  \left\{
     \begin{array}{l}
       L\subset K \\
       \mbox{(subfield)}
     \end{array}
     \;\vert\;
     \begin{array}{l}
       K/L \mbox{ is purely inseparable of exponent} \leq 1, \\
       (K:L) < \infty
     \end{array}
     \right\}\\
     \mbox{and} &&\\
  {\cal D} & := &
    \{{\cal L}:\mbox{$p$-Lie algebra of derivations over $K$}\;\vert\;  \dim_K{\cal L}<\infty\}
  \end{eqnarray*}
Namely,
\begin{equation*}
  \begin{array}{rcl}
    {\cal D}  &  \overset{1:1}{\longleftrightarrow} & {\cal E} \\
    {\cal L}      & \To & K^{\cal L}=\{a\in K\;\vert\; D(a)=0
    \; (\forall D\in {\cal L})\}\\
    {\cal L}_L=\{D:K/L\to K/L:\mbox{derivation}\} & \longleftarrow & L \\
    \end{array}
  \end{equation*}

Jacobson's theory has been reformulated to foliation of algebraic varieties
by T.~ Ekedahl~\cite{Ek87}.

\begin{definition}[foliation in characteristic $p$]
  Let $X$ be a normal variety over an algebraically closed field $k$
  of $\chara(k)=p>0$. Then
  a coherent subsheaf
  ${\cal F}\subset T_X$ is called a {\em foliation}
  if it is closed under Lie bracket and $p$th power
\end{definition}

\begin{definition}[quotient by a foliation]
 Let ${\cal F}\subset T_X$ be a foliation. Set
  \begin{equation*}
    Ann({\cal F})(U):= \{a\in \OO_X(U)\;\vert\; D(a)=0\; (\forall D\in
    {\cal F}(U)\}\qquad (\forall U\subset X:\mbox{open})
  \end{equation*}
  Then we call $Y:=\Spec Ann({\cal F})$ the {\em quotient} of $X$
  by the foliation ${\cal F}$.
\end{definition}

  \begin{proposition}
    \label{ek:foliation}
    Under the above hypothesis, we have
    \begin{enumerate}
      \item [$(i)$] if $\rk {\cal F}=r$, then $[k(X):k(Y)]= p^r$,
      \item [$(ii)$] if $X$ is smooth, $Y$ is also smooth if and only
        if ${\cal F}$ is a subbundle, i.e., $T_X/{\cal F}$ is
        locally free.
      \end{enumerate}
    \end{proposition}

\begin{definition}[$\Sing({\cal F})$]
  Based on Prop.~\ref{ek:foliation}(ii), we consider
  the maximal Zarisiki open subset $U\subset X$
  such that ${\cal F}$ is a subbundle of $T_X\vert_U$.
  Then set $\Sing({\cal F}) = X-U$, which we call
  the {\em singular locus} of the foliation ${\cal F}$.
\end{definition}

$\Sing({\cal F})$ is the points of $X$ that will be
the singular locus of the quotient $\Spec(Ann({\cal F}))$.

\subsubsection{construction of Hirokado~99 variety}

We show an outline of the construction of 
a non-liftable Calabi-Yau 3-fold in characteristic~$3$ \cite{H99}.
This is an application of quotient by foliation.

\begin{description}
  \item [Step 1]
Consider the derivation
\begin{equation*}
  \delta := (x^p-x)\frac{\partial}{\partial x}
  +(y^p-y)\frac{\partial}{\partial y}
  +(z^p-z)\frac{\partial}{\partial z}.
\end{equation*}
on the affine
open subset ${\Bbb A}^3:=\Spec k[x,y,z] \subset {\Bbb P}^3$,
which can be naturally extended to $\PP^3$.
We know that $\delta^p = (-1)^{p-1}\delta$, so that we obtain a
rank one foliation ${\cal F}:=\langle\delta\rangle
(\subset T_{\PP^3})$.
We know that  $\Sing({\cal F}) = \PP^3_{\FF_p}(\subset \PP^3)$, which are
    $p^3+p^2+p+1$ isolated points.

\item [Step 2]
  Each singular point in $\Sing({\cal F})$ can
  be resolved by one point blow-up.
  Then we have the following diagram
  \begin{equation*}
  \begin{CD}
    S         @>{g}>>    X               @>{\overline{g}}>>  S^{(p)} \\
    @V{\pi}VV         @V{\overline{\pi}}VV   @VVV\\
      \PP^3 @>{g_0}>> V @>{\overline{g}_0}>> (\PP^3)^{(p)}
    \end{CD}
\end{equation*}where (1) $\pi:S\To {\Bbb P}^3$ is the blow-ups at $p^3+p^2+p+1$ singular
points in $\Sing({\cal F})$,
(2) $g$  is the quotient morphism (of  degree $p$)
induced by the foliation $\pi^*{\cal F}$, which has no singular point,
(3) $g_0$ is the quotient morphism (of degree $p$) induced by
${\cal F}$,
(4) $\overline{\pi}$ is the naturally induced birational morphism, which
is a desingularization of $V$, 
and (5) $\overline{g}$ and $\overline{g}_0$
are the  morphisms of degree $p^2$ in the 
relative Frobenius morphisms  $\overline{g}\circ g$ and
$\overline{g}_0\circ g_0$.

\item [Step 3] We have
\begin{eqnarray*}
  g^*\omega_X
  &\iso &\pi^*\OO_{\PP^3}((p-1)^2-4) \tensor \OO_S(D)\\
  \mbox{where} &&
   D = (3-p)\sum_{i=1}^{p^3+p^2+p+1}E_i
  \quad (E_i:\mbox{exceptional divisor of }\pi).
\end{eqnarray*}
In particular, for $p=3$, we have
$g^*\omega_X\iso \pi^*\OO_{\PP^3}$ from which
we obtain $\omega_X\iso \OO_X$ together with 
$H^1(X,\OO_X)=H^2(X,\OO_X)=0$. Thus $X$ is a Calabi-Yau 3-fold.\\
\end{description}

This example has the following properties:
(1) $b_3(X)=0$ so that non-liftable and supersingular,
(2) $X$ is unirational and hence simply connected,
(3) Hodge symmetry holds
and moreover (4) weak $H^1$-Kodaira vanishing holds, i.e., $H^1(X,L^{-1})=0$
for ample $L$ such that $H^0(X,L^p)\ne{0}$
\cite{T1}, however,
according to \cite{Ek03}, it is not liftable to $W_2$.

\subsection{construction (II): supersingular K3 pencil over $\PP^1$}

Let $k$ be an algebraically closed field of characteristic $p>0$.
Let $f:X\To \PP^1_k$ be a smooth morphism such that
$\omega_X = \OO_X$ and
all geometric fibers $X_{\bar{t}}$ ($t\in\PP^1_k$) are K3 surfaces with
Picard number $\rho(X_{\hat{t}})=22$, i.e., supersingular K3 surfaces.
Such a variety $X$ is necessarily a projective Calabi-Yau 3-fold
(see Prop.~1.1, Prop.~1.7~\cite{Schr04}).

Moreover, we have $b_3(X)=0$ so that non-liftable to characteristic $0$,
$X$ is unirational and hence simply connected and weak $H^1$-Kodaira
vanishing holds \cite{T2}, however,
according to \cite{Ek03}, it is not liftable to $W_2$.

Such a 3-fold $X$ is called
a {\em Schr\"oer variety}.

\subsubsection{construction of Schr\"oer varieties}
The outline of the Schr\"oer's construction is as follows. The idea
is taking quotient of (a modified version of) Moret-Bailly
construction of relative abelian surfaces over $\PP^1$ \cite{MB} by
involution to obtain Kummer surface pencil if $p=3$.
For $p=2$, we consider a generalized Kummer surface
instead of ordinary Kummer surface.

\begin{description}
\item [Step 1] Let $A = E_1\times E_2$, where $E_i$ are
  supersingular elliptic curves. Then we have a subgroup scheme
  \begin{equation*}
    K(L):= \alpha_p\times\alpha_p\subset A
    \qquad\mbox{where }\alpha_p := \Spec k[X]/(X^p)
  \end{equation*}
  Then  we take the product with $\PP^1$:
  \begin{equation*}
  K(L)\times {\PP^1} = \Spec
  \OO_{\PP^1}[x,y]/(x^p,y^p) \subset A\times \PP^1.
  \end{equation*}

\item [Step 2]   We take
  \begin{equation*}
      H =\Spec \OO_{\PP^1}[x,y]/(x^p,y^p, tx-sy) (\subset K(L)\times\PP^1)
  \end{equation*}
  where $[s:t]$ is the projective
  coordinates of $\PP^1$.
 Actually, $H$ is a group scheme of height one
  corresponding to $p$-Lie subalgebras $\OO_{\PP^1}(-1)\subset
  \OO_{\PP^1}\dirsum\OO_{\PP^1} (\subset Lie(A\times\PP^1))$.
  Now take $\tilde{X}:= (A\times\PP^1)/H$, which is a relative
  abelian surface on $\PP^1$.\\

\item [Step 3] Take the fiberwise quotient by
  \begin{itemize}
    \item the involution induced from $(x,y)\mapsto
      (-x,-y)$ on $A$ together with desingularization (if $p=3$), or
    \item an automorphism
      $\sigma:S\to S$ of order $3$
      together with desingularization (if $p=2$)
  \end{itemize}
  we obtain the (generalized) Kummer
  surface pencil $\pi:X\to \PP^1$, which is a Calabi-Yau threefold.
  \end{description}

\subsection{construction (III): Schoen type examples}

This construction is originated from \cite{SchMZ} and applied to
produce non-liftable Calabi-Yau 3-folds by many authors
including Schoen himself \cite{Scho09}.
Basic ideal is as follows:

\begin{description}
\item [Step 1] Take the fiber product
  \begin{equation*}
    W:=Y_1\times_{\PP^1}Y_2
  \end{equation*}
  where $\pi_i:Y_i\to \PP^1$ are (rational quasi-)elliptic surfaces
  with section (to use Weierstrass form).
  
\item [Step 2] Carry out desingularization of $W$:
  Singular points of $W$ come from 
  singular fibers of $Y_i$. In some cases, the singularities
  of $W$ can be resolved with small resolution.

\item [Step 3] The desingularization $\overline{W}$
of $W$ is a Calabi-Yau 3-fold with desired properties.
\end{description}

\subsubsection{examples by Hirokado-Ito-Saito}

In the fiber product $Y_1\times_{\PP^1}Y_2$,
$\pi_i:Y_i\to \PP^1$ are 
(1) a rational quasi-elliptic surface with section $Y_1$
    and a rational elliptic surface with section $Y_2$ in
    characteristic $p=2,3$ in \cite{H01}
    or  (2)
quasi-elliptic rational surfaces $Y_i$ with section
    in characteristic $p=3$ \cite{HISark} and $p=2$
     \cite{HISmanuscripta}.
    
For these examples, we have 
$b_3(X)=0$ (hence non-liftable),
and  have structures of fibration over $\PP^1$ or $\PP^2$.
Also they are unirational, hence simply connected.

\subsubsection{examples by Cynk, van Straten and Sch\"utt
\label{section:cynkschutt}}

We can apply the theory presented in section~\ref{section:cynk} to
Schoen type construction. Then without using the condition $b_3(X)=0$
we can assure non-liftability. But $b_3(X)=0$ for some of the examples
are proved in \cite{Scho09}. Also, since we use small resolution, most
of the obtained examples are not projective but algebraic spaces.

\subsection{construction (IV):
  double cover of $\PP^3$-- examples by Cynk and van Straten}

Based on the theory presented in section~\ref{section:cynk}, we can
     construct the examples in the following way (at least
     conceptually):

\begin{description}
\item [Step 1 (CY3-fold over $\chara=0$)]
  Let $X_0$ be a Calabi-Yau 3-fold $X_0$ over a field $K$ of
  characteristic $0$. 
  
\item [Step 2 (find a good reduction-- preparation)]
  Let
    ${\mathcal X}\to \Spec(R)$
  be a scheme over
  over a suitable $\ZZ$-algebra $R$ with $K = Q(R)$
  such that $X_0$ is its generic fiber:${\mathcal X}_\eta=X_0$.
  
\item [Step 3 (find a CY reduction)]
  The exists a non-empty Zariski open set $V\subset \Spec(R)$
  such that the special fiber 
 $X_p := {\mathcal X}\times_R k(P)$  
  is non-singular  Calabi-Yau
  for $P\in V$.

\item [Ste 4 (find a rigid singular fiber $X_p$)]
  Now we find $P\in \Spec(R)- V$ such that $X_P$ is Calabi-Yau
  with node as singularity.

\item [Step 5 (small resolution)]
  Now take a desingularization $\phi: Y\To X_P$ with small resolution.
  Then $Y$ is the desired example.
  \end{description}

In fact, we first construct ${\mathcal X}$ and then find $X_0$.
${\mathcal X}$ is obtained by double cover of $\PP^3$ ramified at some
particular octic or Clebsch diagonal cubic.

The examples in \ref{section:cynkschutt} are also constructed
with this method.

\subsection{Raynaud-Mukai construction cannot produce CY 3-folds}

The famous counter-example to Kodaira vanishing
given by Raynaud \cite{Ray} has been generalized by
Mukai \cite{Mu} and we can produce 3-folds
which is non-liftable even over $W_2$
(see Corollary~\ref{KNAposchar}).
Then it was hoped that non-liftable Calabi-Yau 3-folds
on which Kodaira vanishing does not hold could be
produced with this method.
However, it turned out that it is impossible \cite{T0,T2}.

Whether there exists a non-liftable Calabi-Yau variety
on which Kodaira vanishing does not hold seems to be
still open.


%
%
\bigFoot

\end{document}